\documentclass{amsproc}
\usepackage{graphicx}
\usepackage{amssymb}
\usepackage{epstopdf}
\usepackage{amsmath,amsthm,amscd,amssymb}
\usepackage{latexsym}
\usepackage[colorlinks,citecolor=red,pagebackref,hypertexnames=false]{hyperref}
\usepackage{geometry}                
\numberwithin{equation}{section}

\theoremstyle{plain}
\newtheorem{theorem}{Theorem}[section]
\newtheorem{lemma}[theorem]{Lemma}
\newtheorem{corollary}[theorem]{Corollary}
\newtheorem{proposition}[theorem]{Proposition}

\theoremstyle{definition}
\newtheorem{definition}[theorem]{Definition}

\newtheorem{case[theorem]}{Case}

\theoremstyle{remark}
\newtheorem{remark}[theorem]{Remark}

\numberwithin{equation}{section}

\begin{document}

\title{The Fuglede Conjecture holds in ${\Bbb Z}_p \times {\Bbb Z}_p$} 


\author{Alex Iosevich, Azita Mayeli and Jonathan Pakianathan}

\date{today}

\address{Department of Mathematics, University of Rochester, Rochester, NY}
\email{iosevich@math.rochester.edu}
\address{Department of Mathematics, CUNY Queensborough, NY} 
\email{AMayeli@qcc.cuny.edu}
\address{Department of Mathematics, University of Rochester, Rochester, NY}
\email{jonpak@math.rochester.edu}


\maketitle

\begin{abstract} In this paper we study subsets $E$ of ${\Bbb Z}_p^d$ such that any function $f: E \to {\Bbb C}$ can be written as a linear combination of characters orthogonal with respect to $E$. We shall refer to such sets as spectral. In this context, we prove the Fuglede Conjecture in ${\Bbb Z}_p^2$ which says that $E \subset {\Bbb Z}_p^2$ is spectral if and only if $E$ tiles ${\Bbb Z}_p^2$ by translation. Arithmetic properties of the finite field Fourier transform, elementary Galois theory and combinatorial geometric properties of direction sets play the key role in the proof.  \end{abstract} 


\section{Introduction}

\vskip.125in

Let $E \subset {\Bbb Z}_p^d$, where ${\Bbb Z}_p$, $p$ prime, is the cyclic group of size $p$ and ${\Bbb Z}_p$ is the $d$-dimensional vector space over ${\Bbb Z}_p$. We say that $L^2(E)$ has an orthogonal basis of exponentials (indexed by $A$) if the following conditions hold: 

\vskip.125in 

\begin{itemize} 

\item (Completeness) There exists $A \subset {\Bbb Z}_p^d$ such that for every function $f: E \to {\Bbb C}$ there exist complex numbers 
${ \{c_a\} }_{a \in A}$, $A \subset {\Bbb Z}_p^d$, such that 
$$ f(x)=\sum_{a \in A} \chi(x \cdot a) c_a,$$ for all $x \in E$, where $\chi(u)=e^{\frac{2 \pi i}{p}u}$. We shall refer to $A$ as a {\it spectrum} of $E$. 
Note we may apply this expansion to functions $f: \mathbb{Z}_p^d \to \mathbb{C}$ by first restricting them to $E$ but the resulting equality only holds 
for $x \in E$.

\vskip.125in

\item (Orthogonality) The relation 
$$ \widehat{E}(a-a') \equiv p^{-d} \sum_{x \in E} \chi(x \cdot (a-a'))=0$$ holds for every $a,a' \in A$, $a \not=a'$. 

\vskip.125in 

If these conditions hold, we refer to $E \subset {\Bbb Z}_p^d$ as a {\it spectral} set. 

\end{itemize} 

\vskip.125in 

\begin{definition} (Spectral Pair) A {\bf spectral pair} $(E,A)$ in $V=\mathbb{Z}_p^d$ is a spectral set $E$ with an orthogonal basis of exponentials indexed by $A$. \end{definition} 

\begin{definition} (Tiling Pair) A {\bf tiling pair} $(E',A')$ consists of $E',A' \subset {\Bbb Z}_p^d$ such that every element $v \in V$ can be written uniquely as a sum $v=e'+a', e' \in E', t \in A'$. Equivalently, $(E',A')$ is a tiling pair if $\sum_{a' \in A'} E'(x-a') \equiv 1$ for every $x \in V$. We say that $E'$ {\bf tiles} $V$ by translation if there exists $A' \subset V$ such that $(E',A')$ is a tiling pair. Here and throughout $E(x)$ is the indicator function of $E$. \end{definition}

\vskip.125in 

The study of the relationship between exponential bases and tiling has its roots in the celebrated Fuglede Conjecture in ${\Bbb R}^d$, which says that if $E \subset {\Bbb R}^d$ of positive Lebesgue measure, then $L^2(E)$ possesses an orthogonal basis of exponentials if and only if $E$ tiles ${\Bbb R}^d$ by translation. Fuglede proved this conjecture in the celebrated 1974 paper (\cite{Fu74}) in the case when either the tiling set or the spectrum is a lattice. A variety of results were proved establishing connections between tiling and orthogonal exponential bases. See, for example, \cite{LRW00}, \cite{IP98}, \cite{L02}, \cite{KL03} and \cite{KL04}. In 2001, Izabella Laba proved the Fuglede conjecture for unions of two intervals in the plane (\cite{L01}). In 2003, Katz, Tao  and the first listed author (\cite{IKT03}) proved that the Fuglede conjecture holds for convex planar domains. 

A cataclysmic event in the history of this problem took place in 2004 when Terry Tao (\cite{T04}) disproved the Fuglede Conjecture by exhibiting a spectral set in ${\Bbb R}^{12}$ which does not tile. The first step in his argument is the construction of a spectral subset of ${\Bbb Z}_3^5$ of size $6$. It is easy to see that this set does not tile $\mathbb{Z}_3^5$ because $6$ does not divide $3^5$. As a by-product, this shows that spectral sets in ${\Bbb Z}_p^d$ do not necessarily tile. See \cite{KM06}, where Kolountzakis and Matolcsi also disprove the reverse implication of the Fuglede Conjecture. The general feeling in the field was that sooner or later the counter-examples of both implication will cover all dimensions. However, we see in this paper that the Fuglede Conjecture holds in two-dimensional vector spaces over prime fields. 

\vskip.125in 

Our main result is the following. 

\begin{theorem} \label{zeroesgalore} Let $E$ be a subset of $\mathbb{Z}_p^d$, $p$ an odd prime.  \begin{itemize} 

\vskip.125in 

\item i) (Density) The space $L^2(E)$ has an orthogonal basis of exponentials indexed by $A$ if and only if $|E|=|A|$ and $\hat{E}(a-a')=0$ for all distinct $a, a' \in A$. 

\vskip.125in 

\item ii) If  $E \subset {\Bbb Z}_p^d$ is spectral and $|E| > p^{d-1}$ then $E=A=\mathbb{Z}_p^d$.

\vskip.125in 

\item iii) (Divisibility) If $E \subset {\Bbb Z}_p^d$ is spectral, then $|E|$ is $1$ or a multiple of $p$. 

\vskip.125in 

\item iv) (Fuglede Conjecture in ${\Bbb Z}_p^2$) A set $E \subset {\Bbb Z}_p^2$ is a spectral set if and only if $E$ tiles ${\Bbb Z}_p^2$ by translation. 

\end{itemize} 

\end{theorem} 

\vskip.125in 

\begin{remark} It is not difficult to see that the Fuglede conjecture also holds in ${\Bbb Z}_p$, the one-dimensional setting. This is because a tiling set $E$ has to have order dividing $p$ and hence must be either a point or the whole space which are both trivially spectral sets also. Conversely, a spectral set $E$, by part (iii) of Theorem~\ref{zeroesgalore} has to have order $1$ or a multiple of $p$ and hence must be a point or the whole space which makes it trivially a tiling set also. \end{remark} 

\begin{remark} It is interesting to note that Theorem~\ref{zeroesgalore} also holds for $p=2$. Parts (i)-(iii) are immediate following the same proofs given in this paper for the odd prime case. Part (iv) follows trivially since parts (i)-(iii) imply that the only tiling or spectral sets in 
$\mathbb{Z}_2^2$ have orders $1, 2$ or $4$. In the only non-trivial case, when $|E|=2$, $E$ has to be a line which is easily seen to be both a tiling and a spectral set. \end{remark} 

\vskip.125in 

\subsection{Acknowledgements} The authors wish to thank Mihalis Kolountzakis and Michael Lacey for several helpful remarks. 

\vskip.25in 

\section{Basic properties of spectra} 

\vskip.125in 

\begin{lemma} \label{coefficients} Suppose that $L^2(E)$ has an orthogonal basis of exponentials and 
$$f: {\Bbb Z}_p^d \to {\Bbb C}.$$ Then the coefficients are given by 
$$ c_a(f)={|E|}^{-1} \sum_{x \in E} \chi(-x \cdot a) f(x).$$ 

\end{lemma} 

To prove this, observe that if $f(x)=\sum_{a \in A} \chi(x \cdot a) c_a$ for $x \in E$, then 

\begin{eqnarray*}
{|E|}^{-1} \sum_{x \in E}  \chi(-x \cdot a)f(x) &=& {|E|}^{-1} \sum_{x \in E} \sum_{b \in A} \chi(-x \cdot (a-b))c_b(f) \\
&=& {|E|}^{-1} \sum_{b \in A} c_b(f) \sum_{x \in E} \chi(-x \cdot (a-b))=c_a(f)
\end{eqnarray*}

 and the proof is complete. 

\vskip.125in 

\begin{lemma} (Delta function) \label{delta} Suppose that $L^2(E)$ has an orthogonal basis of exponentials with the spectrum $A$. Let 
$\delta_0(x)=1$ if $x=\vec{0}$ and $0$ otherwise. Furthermore assume $\vec{0} \in E$. Then 
$$ \delta_0(x)={|E|}^{-1} \sum_{a \in A} \chi(x \cdot a).$$ 

\end{lemma} 

To prove the lemma, observe that if $f(x)=\delta_0(x)$, then 
$$ c_a(f)={|E|}^{-1} \sum_{x \in E} \chi(-x \cdot a) \delta_0(x)={|E|}^{-1}.$$ 

The conclusion follows from Lemma \ref{coefficients}. 

\vskip.125in 

\begin{lemma} (Parseval) \label{parseval} Suppose that $L^2(E)$ has an orthogonal basis of exponentials and $f$ is any function on 
${\Bbb Z}_p^d$ with values in ${\Bbb C}.$ Then 
$$ \sum_{a \in A} {|c_a(f)|}^2={|E|}^{-1} \sum_{x \in E} {|f(x)|}^2.$$ 
\end{lemma} 

\vskip.125in 

\begin{lemma} (Density) \label{density} Suppose that $L^2(E)$ has an orthogonal basis of exponentials with the spectrum $A$. 
Then $|E|=|A|$. 
\end{lemma} 

\vskip.125in 

To see this, observe that since $\{ \chi(x \cdot a): a \in A \}$ is an orthogonal set of nonzero elements for $L^2(E)$ we know it is a linearly independent set in $L^2(E)$. By completeness we know it also spans $L^2(E)$ and hence is a basis for $L^2(E)$. Thus 
$|A| = \dim_{\mathbb{C}}(L^2(E)) = |E|$.


\vskip.25in 

\section{Proof of Theorem \ref{zeroesgalore}} 

\vskip.125in 

To prove Part i) of Theorem \ref{zeroesgalore} we note that the orthogonality property is equivalent to $\widehat{E}(a-a')=0$ for all $a \not=a'$, $a,a' \in A$. Recall that since $\{ \chi(x \cdot a): a \in A \}$ is linearly independent in $L^2(E)$, it spans $L^2(E)$ if and only if 
$|E|=|A|$.

\begin{definition} (\cite{IMP12}) We say that two vectors $x$ and $x'$ in ${\Bbb Z}_p^d$ point in the same direction if there exists $t \in {\Bbb Z}_p^{*}$ such that $x'=tx$. Here ${\Bbb Z}_p^{*}$ denotes the multiplicative group of ${\Bbb Z}_p$. Writing this equivalence as $x \sim x'$, we define the set of directions as the quotient 
\begin{equation} \label{directionequivalence} {\mathcal D}({\Bbb Z}_p^d)={\Bbb Z}_p^d / \sim. \end{equation}

Similarly, we can define the set of directions determined by $E \subset {\Bbb Z}_p^d$ by 

\begin{equation} \label{equivalence} {\mathcal D}(E)=E-E / \sim, \end{equation} where 
$$ E-E=\{x-y: x,y \in E\},$$ with the same equivalence relation $\sim$ as in (\ref{directionequivalence}) above. 
\end{definition} 

\vskip.125in 

\begin{theorem} \label{maindirection} (\cite{IMP12}) Let $E \subset {\Bbb Z}_p^d$.  Suppose that $|E|>p^k$, $1 \leq k \leq d-1$. Let $H_{k+1}$ denote a $k+1$-dimensional sub-space of ${\Bbb Z}_p^d$. Then ${\mathcal D}(H_{k+1}) \subset {\mathcal D}(E)$. In particular, if $|E|>p^{d-1}$, every possible direction is determined. 

Furthermore a set $E$ does not determine all directions if and only if there is a hyperplane $H$ and $S \subseteq H$ such that $E$ is the graph of a function $f: S \to {\Bbb Z}_p$ over $H$, which means that relative to some decomposition ${\Bbb Z}_p^d = H \oplus {\Bbb Z}_p$, $E=\{ (x,f(x)): x \in S \}$.

\end{theorem}

\vskip.125in 

\begin{theorem}(\cite{HIPRV15})
\label{huy} 
Let $E \subset {\Bbb Z}_p^d$. Then $\hat{E}(m) = 0$ implies that $\hat{E}(rm)=0$ for all 
$r \in \mathbb{Z}_p^{*}$. Furthermore $\hat{E}(m)=0$ for $m \neq 0$ if and only if $E$ is equidistributed on the $p$ hyperplanes 
$H_t=\{x: x \cdot m = t \}$ for $t \in {\Bbb Z}_p$ in the sense that 
$$ \sum_{x \cdot m=t} E(x)=|E \cap H_t|,$$ viewed as a function of $t$, is constant. \end{theorem}

\begin{remark} Note this last theorem is a fact about rational valued functions over prime fields that is {\bf not true} for complex valued functions in general or over other fields. We give the proof of Theorem \ref{huy} at the end of the paper for the sake of completeness. \end{remark} 

\vskip.125in 

The proof of part ii) of Theorem \ref{zeroesgalore} follows fairly easily from combining Theorem \ref{maindirection} and Theorem \ref{huy}. Indeed, suppose that $L^2(E)$ has an orthogonal basis of exponentials and $|E|>p^{d-1}$. By Lemma \ref{density}, $|E|=|A|>p^{d-1}$. By Theorem \ref{maindirection}, ${\mathcal D}(A)={\mathcal D}({\Bbb Z}_p^d)$. Combining this with Theorem \ref{huy} implies that $\widehat{E}$ vanishes on ${\Bbb Z}_p^d \ \backslash \vec{0}$. It follows that $E={\Bbb Z}_p^d$, as claimed. To see this, observe that if we define for $f: {\Bbb Z}_p^d \to {\mathbb C}$, 
\begin{equation} \label{ftdef} \widehat{f}(m)=p^{-d} \sum_{x \in {\Bbb Z}_p^d} \chi(-x \cdot m) f(x), \end{equation} then 
\begin{equation} \label{ftinvdef} f(x)=\sum_{m \in {\Bbb Z}_p^d} \chi(x \cdot m) \widehat{f}(m). \end{equation} It follows that $E(x)=\widehat{E}(0, \dots, 0)=|E|p^{-d}$. Since $E(x)$ is an indicator function of a set, we conclude that $|E|=p^d$ and the claim is proved. 

\vskip.125in 

We define a trivial spectral pair in $\mathbb{Z}_p^d$ to be $(E,A)=(\text{point}, \text{another point})$ or $(E,A)=(\mathbb{Z}_p^d, \mathbb{Z}_p^d)$ or $(E,A) = (\emptyset, \emptyset)$.

\vskip.125in 

Part iii) of Theorem \ref{zeroesgalore} is contained in the following result. 

\begin{proposition} \label{propdivide} 
\label{theorem: Fugledeprep1} Let $p$ be an odd prime and $(E,A)$ be a non-trivial spectral pair in $\mathbb{Z}_p^d$ then $|E|=|A|=mp$ where $m \in \{1,2,3,\dots, p^{d-2} \}$. 
\end{proposition}

To prove Proposition \ref{theorem: Fugledeprep1}, let $(E,A)$ be a non-trivial spectral pair. Then part i) of 
Theorem~\ref{zeroesgalore} shows that 
$|E|=|A|$ and $\hat{E}(a-a')=0$ for distinct $a, a' \in A$. Since the spectral pair $(E,A)$ is nontrivial, $2 \leq |E|=|A| \leq p^{d-1}$ also. Thus taking two distinct elements $a, a' \in A$ shows that $\hat{E}(\alpha)=0$ for 
$\alpha = a -a' \neq 0$. Thus $E$ is equidistributed on the $p$ parallel hyperplanes 
$$H_t = \{x: x \cdot \alpha = t \},$$ $t \in \mathbb{Z}_p$ by Theorem~\ref{huy}. Thus if $E$ has $m \geq 1$ elements per hyperplane we have $|E|=|A|=mp$. Then $1 \leq m \leq p^{d-2}$ since $0 < |E| \leq p^{d-1}$. This proves part iii) of Theorem \ref{zeroesgalore}. 

\vskip.125in 

Observe that if $d=2$ and $E$ is a (non-trivial) spectral set, then $|E|=|A|=mp$ implies $|E| \geq p$ while $|E| \leq p$ by part ii) of Theorem~\ref{zeroesgalore} and so $|E|=|A|=p$. Futhermore, by Theorem \ref{maindirection} above, $A$ is a graph of a function ${\Bbb Z}_p \to {\Bbb Z}_p$ since $|A|=p$ and it does not determine all directions. Finally, since $E$ is equidistributed on a family of $p$ parallel lines and $|E|=p$, we see that $E$ is also a graph of a function ${\Bbb Z}_p \to {\Bbb Z}_p$ with respect to some system of axes. The following is an immediate corollary of Proposition \ref{propdivide}. 

\begin{corollary} If $E$ is a spectral set in $\mathbb{Z}_p^2$, $p$ an odd prime, then $E$ is either a point, a graph set of order 
$p$ or the whole space and hence tiles $\mathbb{Z}_p^2$ in all cases.
\end{corollary}

This corollary follows from Proposition~\ref{theorem: Fugledeprep1} immediately once one notes that any graph set 
$$E=\{(x, f(x)): x \in \mathbb{Z}_p \}$$ for a function $f$, with respect to some coordinate systems, tiles $\mathbb{Z}_p^2$ using the tiling partner 
$$A = \{ (0,t): t \in \mathbb{Z}_p \}.$$

To complete the proof of the Fuglede conjecture in two dimensions over prime fields, which is the content of part iv) of Theorem \ref{zeroesgalore}, it remains to show that any tiling set is spectral since we have just shown that any spectral set tiles. 
 
\begin{proposition} (Sets which tile by translation are spectral) 
Let $p$ be an odd prime, and let $E \subseteq \mathbb{Z}_p^2$. Suppose that $E$ tiles $\mathbb{Z}_p^2$ by translation. Then $E$ is a spectral set. \end{proposition}

We shall need the following result from  \cite{HIPRV15}. We shall prove it at the end of the paper for the sake of completeness. 

\begin{theorem} \label{tileshot} (\cite{HIPRV15}) Let $E$ be a set that tiles $\mathbb{Z}_p^2$. Then $|E|=1,p$ or $p^2$ and $E$ is a graph if $|E|=p$. \end{theorem} 

The cases $|E|=1, p^2$ are trivially spectral sets so we may reduce to the case that $E$ is a graph, i.e
$$E=\{ xe_1 + f(x)e_2: x \in \mathbb{Z}_p \}$$ where $e_1, e_2$ is a basis for $\mathbb{Z}_p^2$ and $f: \mathbb{Z}_p \to \mathbb{Z}_p$ is a function. By changing the function if necessary we can assume $e_2$ is orthogonal to $e_1$ unless $e_1 \cdot e_1 = 0$ i.e., unless $e_1$ generates an isotropic line. 
This case does not occur if $p \equiv 3 \mod(4)$. In the case when $p \equiv 1 \mod(4)$, it is possible that $e_1$ generates one of two isotropic lines  
$$ \{(t,it): t \in {\Bbb Z}_p \},$$ where $i$ is one of the two distinct solution of the equation $x^2+1=0$. The reason this case needs to be treated separately is that $(t_1, it_1) \cdot (t_2, it_2)=0$ for all $t_1,t_2 \in {\Bbb Z}_p$. To deal with this, we note that the other solution of the equation $x^2+1=0$ is given by $-i$ and we take $e_2$ on the other isotropic line in the plane, given by 
$$ \{(t,-it): t \in {\Bbb Z}_p \}.$$ Then $e_1 \cdot e_1 = 0 = e_2 \cdot e_2$ and we may assume $e_1 \cdot e_2 = 1$ by scaling $e_2$ appropriately.

\vskip.125in 

There are two situations to consider. 

\begin{itemize} 

\item Case 1: $e_1$ and $e_2$ are orthogonal. Then we will take $A=\{ xe_1: x \in \mathbb{Z}_p \}$. To show that $(E,A)$ is a spectral pair, we need only show $\{ \chi(ae_1 \cdot x): a \in \mathbb{Z}_p \}$ are orthogonal in $L^2(E)$. By Theorem \ref{huy} this happens if and only if $\hat{E}((a-a')e_1) = 0$ for all distinct $a, a' \in \mathbb{Z}_p$ which happens if and only if $E$ equidistributes on the $p$ parallel lines normal to $e_1$, i.e., on the $p$ parallel line of constant $e_1$-coordinate in the $(e_1,e_2)$-grid. This is clearly the case as $E$ is a graph over the $e_1$ coordinate and hence has exactly $1$ element on each of these parallel lines, so this case is proven.

\vskip.125in 

\item Case 2: $e_1$ and $e_2$ generate the two isotropic lines in $\mathbb{Z}_p^2$, $p=1$ mod $4$. In this case $e_1 \cdot e_2 \neq 0$ but $e_1 \cdot e_1 = e_2 \cdot e_2 = 0$. Since $E$ is equidistributed along the $p$ parallel lines of constant $e_1$-coordinate, it is easy to see that these are the same family of lines as $H_t = \{ x: x \cdot e_2 = t \}$, $t \in \mathbb{Z}_p$.  Thus in this case using 
$A= \{ a e_2: a \in \mathbb{Z}_p \}$ we find that $\hat{E}((a-a')e_2) = 0$ for distinct $a, a' \in \mathbb{Z}_p$ and so $(E,A)$ is a spectral pair. Thus $E$ is still spectral in this case and the theorem is proven in all cases.
\end{itemize} 

\vskip.25in 

\section{Proof of Theorem \ref{huy}} 

\vskip.125in 

We have 
$$ \widehat{E}(m)=p^{-d} \sum_{x \in {\Bbb Z}_p^d} \chi(-x \cdot m) E(x)=0$$ for some $m \not=(0, \dots, 0)$. Let $\xi=\chi(-1)$. It follows that 
$$ 0=\sum_{x \in {\Bbb Z}_p^d} \xi^{x \cdot m} E(x)=\sum_{t \in {\Bbb Z}_p} \xi^t \sum_{x \cdot m=t} E(x).$$ 

Let 
$$ n(t)=\sum_{x \cdot m=t} E(x) \in \mathbb{Q},$$ so 
$$ \sum_{t \in {\Bbb Z}_p} \xi^t n(t)=0.$$ 

This means that $\xi$ is a root of the rational polynomial 
$$ P(u)=\sum_{t \in {\Bbb Z}_p} n(t) u^t.$$ 

The minimal polynomial over $\mathbb{Q}$ of $\xi$ is 
$$ Q(u)=1+u+\dots+u^{p-1},$$ so by elementary Galois theory, $P(u)$ is a constant multiple of $Q(u)$ since $P$ is a rational polynomial with $\xi$ as root, $Q$ is the minimal polynomial of $\xi$ and $P$ and $Q$ are both of degree $p-1$. It follows that the coefficients of $n(t)$ are independent of $t$. This proves the second assertion of Theorem \ref{huy}, namely that $E$ is equidistributed on the hyperplanes $H_t=\{x \in {\Bbb Z}_p^d: x \cdot m=t \}$. 

Let us now prove that if $\widehat{E}(m)=0$ for some $m \not=(0, \dots, 0)$, then $\widehat{E}(rm)=0$ for all $r \not=0$. We have 
$$ \sum_{x \in {\Bbb Z}_p^d} \chi(-x \cdot rm) E(x)=\sum_{t \in {\Bbb Z}_p} \xi^t \sum_{x \cdot rm=t} E(x)$$
$$=\sum_{t \in {\Bbb Z}_p} \xi^t \sum_{x \cdot m=tr^{-1}} E(x)=\sum_{t \in {\Bbb Z}_p} \xi^t n(r^{-1}t).$$ 

\vskip.125in 

It follows from above that for a fixed $r$, $n(r^{-1}t)$ is independent of $t$. Therefore 
$$ \sum_{t \in {\Bbb Z}_p} \xi^t n(r^{-1}t)=\sum_{t \in {\Bbb Z}_p} \xi^t n(t)=0$$ and the proof of the claim follows. This completes the proof of Theorem \ref{huy}. 

\vskip.125in 

\begin{remark} We again note that the last theorem can be extended to a theorem about rational valued functions but not in general to one about real valued or complex valued functions. For example the polynomial $P$ was a multiple of the minimal polynomial $Q$ only because it was a rational polynomial. Indeed $\xi$ is the root of $x-\xi$ over $\mathbb{C}$ and $(x-\xi)(x-\bar{\xi}) = x^2 - 2cos(2\pi/p) + 1$ over $\mathbb{R}$ which are not multiples of $Q$ when $p > 3$. Indeed, Galois theory can be used to say more about the relationship between the Fourier coefficients of a rational function along a line, aside from the origin, they are Galois conjugates over some cyclotomic field) and these relationships do not hold for complex valued functions nor over finite fields other than prime fields. See \cite{HIPRV15} for more details. The structure of tiling sets and spectral sets established over prime fields in this paper fundamentally owes  to this fact and the lack of this structure over other fields like $\mathbb{R}$ is part of the reason why the corresponding questions remain unresolved. 
\end{remark} 

\vskip.125in 

\section{Proof of Theorem \ref{tileshot}} 

\vskip.125in 

Let $T$ denote the set that tiles $E$. Note that $|E||T|=p^2$, so $|E|=1, p$ or $p^2$. If $|E|=1$ then $E$ is a point and we are done. If $|E|=p^2$ then $E$ is the whole plane and we are done, so without loss of generality $|E|=p$.

If $\hat{E}(m)$ never vanishes then $E$ is a point and we are done. To see this, observe that by definition, tiling means that 
$$ 1 \equiv \sum_{\tau \in T} E(x-t)=\sum_{t \in {\Bbb Z}_p^d} E(x-\tau)T(\tau)=E \star T(x)$$ for all $x \in {\Bbb Z}_p^d$. 

This immediately implies $\hat{T}(m)\hat{E}(m)=0$ for all nonzero $m$ as the Fourier transform of the convolution is the product of the Fourier transforms. 
Thus $\hat{T}(m)=0$ for all nonzero $m$ as $\hat{E}(m)$ does not vanish.
This in turn means $|T|=p^2$ and so $|E|=1$ as $|E||T|=p^2$.


On the other hand if 
$\hat{E}(m)=0$ for some nonzero $m$, then it vanishes on $L$, the line passing through the origin and $m \not=\vec{0}$
.
Thus if we set $L^{\perp}$ to be the line through the origin, perpendicular to $m$, we see that 
$$\widehat{L^{\perp}}(s)\widehat{E}(s)=0$$ for all nonzero $s$. This is because by a straightforward calculation 
$$\widehat{L^{\perp}}(s)=p^{-(d-1)} L(s).$$ Since $|L^{\perp}|=p=|E|$ we then see that $E$ tiles ${\Bbb F}_p^2$ by $L^{\perp}$. 

Since $\widehat{E}(m)=0$ for some non-zero vector $m$, we see that $E$ is equidistributed on the set of $p$ lines 
$H_t=\{x: x \cdot m = t\}, t \in \mathbb{F}_p$. Since $|E|=p$ this means there is exactly one point of $E$ on each of these lines. 

We will now choose a coordinate system in which $E$ will be represented as a graph of a function. The coordinate system will either be an orthogonal system or an isotropic system depending on the nature of the vector $m$. There are two cases to consider. \begin{itemize} 

\vskip.125in 

\item Case 1: $m \cdot m \neq 0$: We may set $e_1 = m$ and $e_2$ a vector orthogonal to $m$. $\{ e_1, e_2 \}$ is an orthogonal basis because $e_2$ does not lie on the line throough $m$ as this line is not isotropic. If we take a general vector $\hat{x} = x_1 e_1 + x_2 e_2$ we see that $\hat{x} \cdot m = x_1 (m \cdot m)$ 
and so the lines $H_t, t \in \mathbb{F}_p$ are the same as the lines of constant $x_1$ coordinate with respect to this orthogonal basis $\{ e_1, e_2 \}$. Thus there is a unique value of $x_2$ for any given value of $x_1$ so that $x_1 e_1 + x_2 e_2 \in E$. Thus 
$E=\{ x_1 e_1 + f(x_1) e_2: x_1 \in \mathbb{Z}_p \}=Graph(f)$ for some function $f: \mathbb{F}_p \to \mathbb{Z}_p$.

\vskip.125in 

\item Case 2: $m \cdot m = 0$: We may set $e_1 = m$. In this case any vector orthogonal to $e_1$ lies on the line generated by $e_1$ and so cannot be part of a basis 
with $e_1$. Instead we select $e_2$ off the line generated by $e_1$ and scale it so that $e_1 \cdot e_2 = 1$. Then by subtracting a suitable multiple of $e_1$ from $e_2$ one can also ensure $e_2 \cdot e_2 = 0$. Thus $\{ e_1, e_2 \}$ is a basis 
consisting of two linearly independent isotropic vectors. With respect to this basis, the dot product is represented by the matrix 
$$ \mathbb{A} = \begin{bmatrix} 0 & 1 \\ 1 & 0 \end{bmatrix}$$ which exhibits the plane as the hyperbolic plane. This case can only occur when $p=1$ mod $4$. \end{itemize} 

Note when we express a general vector $x=x_1 e_1 + x_2 e_2$ with respect to this basis we have $x \cdot m = x_2$ thus the lines  
$\{H_t: t \in \mathbb{Z}_p \}$ are the same as the lines of constant $x_2$ coordinate with respect to this basis and $E$ has a unique point on each of these lines. Thus 
$$E=\{ f(x_2) e_1 + x_2 e_2: x_2 \in \mathbb{Z}_p \} = Graph(f)$$ is a graph with respect to this isotropic coordinate system.

Finally we note any function $f: \mathbb{Z}_p \to \mathbb{F}_p$ is given by a polynomial of degree at most $p-1$, explicitly expressed in the form 
$$ f(x) = \sum_{k \in \mathbb{Z}_p} f(k) \frac{\Pi_{j \neq k} (x-j)}{\Pi_{j \neq k} (k-j)}$$
expresses $f$ as such a polynomial in $x$.

\vskip.125in 

\end{document}